\documentclass[12pt]{article}
\usepackage{epic}
\begin{document}
%
\newfont{\twelvemsb}{msbm10 scaled\magstep1}
\newfont{\eightmsb}{msbm8} \newfont{\sixmsb}{msbm6} \newfam\msbfam
\textfont\msbfam=\twelvemsb \scriptfont\msbfam=\eightmsb
\scriptscriptfont\msbfam=\sixmsb \catcode`\@=11
\def\Bbb{\ifmmode\let\next\Bbb@\else \def\next{\errmessage{Use
      \string\Bbb\space only in math mode}}\fi\next}
\def\Bbb@#1{{\Bbb@@{#1}}} \def\Bbb@@#1{\fam\msbfam#1}
\newfont{\twelvegoth}{eufm10 scaled\magstep1}
\newfont{\tengoth}{eufm10} \newfont{\eightgoth}{eufm8}
\newfont{\sixgoth}{eufm6} \newfam\gothfam
\textfont\gothfam=\twelvegoth \scriptfont\gothfam=\eightgoth
\scriptscriptfont\gothfam=\sixgoth \def\frak{\frak@}
\def\frak@#1{{\fam\gothfam{{#1}}}} \def\frak@@#1{\fam\gothfam#1}
\catcode`@=12
%
\def\CC{{\Bbb C}}
\def\NN{{\Bbb N}}
\def\QQ{{\Bbb Q}}
\def\RR{{\Bbb R}}
\def\ZZ{{\Bbb Z}}
\def\cA{{\cal A}}          \def\cB{{\cal B}}          \def\cC{{\cal C}}
\def\cD{{\cal D}}          \def\cE{{\cal E}}          \def\cF{{\cal F}}
\def\cG{{\cal G}}          \def\cH{{\cal H}}          \def\cI{{\cal I}}
\def\cJ{{\cal J}}          \def\cK{{\cal K}}          \def\cL{{\cal L}} 
\def\cM{{\cal M}}          \def\cN{{\cal N}}          \def\cO{{\cal O}}
\def\cP{{\cal P}}          \def\cQ{{\cal Q}}          \def\cR{{\cal R}} 
\def\cS{{\cal S}}          \def\cT{{\cal T}}          \def\cU{{\cal U}}
\def\cV{{\cal V}}          \def\cW{{\cal W}}          \def\cX{{\cal X}}
\def\cY{{\cal Y}}          \def\cZ{{\cal Z}}
\def\csl#1{C_{\uq{#1}}}
\def\expo#1{\mbox{\footnotesize #1}}
\def\fracdis#1#2{\frac{\displaystyle #1}{\displaystyle #2}}
\def\GZ{G.--Z. }
\def\Ket#1#2#3#4#5#6{\left| \begin{array}{ccccc} #1 && #2 && #3 \\
    & #4 && #5 & \\ && #6 &&  \end{array} \right.\rrangle}
\def\ket#1#2#3{\left| \begin{array}{ccc} 
     #1 && #2  \\ & 
     \hspace {-.5cm} #3 \hspace {-.5cm} &  \end{array} \right\rangle}
\def\keti{\ket{p_{12}}{p_{22}}{p_{11}}}
\def\qed{\hfill \rule{5pt}{5pt}}
\def\fdirreps{\def\fdirreps{f.~d. irreps}finite dimensional irreducible
  representations (f.~d. irreps)}
\def\qmbox#1{\qquad\mbox{#1}\qquad}
\def\rrangle{
  \setlength{\unitlength}{0.01cm}
  \begin{picture}(20,120)(0,0)
    \put(0,-62){\line(1,4){18}}
    \put(0,82){\line(1,-4){18}}
  \end{picture}}
\def\uq#1{{\cal U}_q(sl(#1))}
\newtheorem{lemma}{Lemma}
\newtheorem{prop}{Proposition}
\newtheorem{theo}{Theorem}
\newenvironment{result}{\vspace{.2cm} \em}{\vspace{.2cm}}
\renewcommand{\thefootnote}{\fnsymbol{footnote}}
\newpage
\pagestyle{empty}
\setcounter{page}{0}
%
\newcommand{\norm}[1]{{\protect\normalsize{#1}}}
\newcommand{\LAP}{{\small E}\norm{N}{\large S}{\Large L}%
  {\large A}\norm{P}{\small P}}
\newcommand{\sLAP}{{\scriptsize E}{\footnotesize{N}}{\small S}%
  {\norm L}{\small A}{\footnotesize{P}}{\scriptsize P}}
\def\logolapin{
  \raisebox{-1.2cm}{\epsfbox
    {/lapphp8/keklapp/ragoucy/paper/enslapp.ps}}}
\def\logolight{{\bf {\large E}{\Large N}{\LARGE S}{\huge L}%
    {\LARGE A}{\Large P}{\large P} }}
\def\logoenslapp{\logolight}
%
%
%
\hbox to \hsize{
\hss
\begin{minipage}{5.2cm}
  \begin{center}
    {\bf Groupe d'Annecy\\ \ \\
      Laboratoire d'Annecy-le-Vieux de Physique des Particules}
  \end{center}
\end{minipage}
\hfill
\logoenslapp
\hfill
\begin{minipage}{4.2cm}
  \begin{center}
    {\bf Groupe de Lyon\\ \ \\
      {\'E}cole Normale Sup{\'e}rieure de Lyon}
  \end{center}
\end{minipage}
\hss}

\vspace {.3cm}
\centerline{\rule{12cm}{.42mm}}
%
\vspace{20mm}

\begin{center}

  {\LARGE  {\bfseries\sffamily 
      Non integrable  
      representations of the restricted 
      quantum analogue of $sl(3)$ at roots of 1
      }}\\[1cm]
  
  \vspace{10mm}
  
  {\large D. Arnaudon \footnote{arnaudon@lapp.in2p3.fr.
      \\ \indent \ \ Some pictures are available at 
      http://lapphp0.in2p3.fr/\~{}arnaudon/teepee.html
      }}\\[.42cm]
  
  {\em Laboratoire de Physique Th{\'e}orique }\LAP\footnote{URA 14-36 
    du CNRS, associ{\'e}e {\`a} l'{\'E}cole Normale Sup{\'e}rieure de Lyon et
    {\`a} l'Universit{\'e} de Savoie.}, CNRS\\[.242cm]
  
  Groupe d'Annecy: LAPP, BP 110, F-74941
  Annecy-le-Vieux Cedex, France.
  \\

\end{center}
\vspace{20mm}

\begin{abstract}
  The structure of irreducible representations of (restricted)
  $\cU_q(sl(3))$ at roots of unity is understood within the
  Gelfand--Zetlin basis. The latter needs a weakened definition for
  non integrable representations, 
  where the quadratic Casimir operator of the quantum subalgebra 
  $\cU_q(sl(2))\subset\cU_q(sl(3))$ is not completely diagonalized. 
  This is necessary in order to
  take in account the indecomposable $\cU_q(sl(2))$-modules that
  appear. The set of redefined (mixed) states has a teepee shape
  inside the pyramid made with the whole representation.
\end{abstract}

\vfill
\rightline{q-alg/9610025}
\rightline{\LAP-A-625/96}
\rightline{October 1996}

\newpage
\pagestyle{plain}
\renewcommand{\thefootnote}{\arabic{footnote}}

\section{\label{sect:introduction}Introduction}

In this paper, we are interested in 
finite dimensional representations of the quantum
analogue of the enveloping algebra of $sl(3)$ at roots of unity, in
the restricted specialization. 

When the deformation parameter $q$ is not a root of unity, the finite
dimensional irreducible
representations of quantum groups as defined in \cite{Dri,Jim} 
are in correspondence with the
classical ones \cite{RosA,Lusztigi}. 
This correspondence is $2^{\expo{rank}}$-to-one, 
the factor 
$2^{\expo{rank}}$ being
related to trivial isomorphisms of the quantum enveloping algebra.

When $q$ is a root of unity, the dimension of finite dimensional
irreducible representations is bounded. 
In the unrestricted specialization, new classes of
irreducible representations 
appear,  that 
are characterized by continuous parameters (See \cite{RA} for $\uq2$ and
\cite{DK} for general $\cU_q(\cG)$). 
We do not consider them here, since we are interested in the restricted
specialization, and more precisely in its finite dimensional Hopf
subalgebra, where the raising and lowering generators are
nilpotent and where the Cartan generators are quantized. 
In this case, the finite dimensional irreducible  
representations can be obtained as quotient of Verma modules 
(with integral dominant highest weights)
by their maximal submodule. 

\medskip

As for representations of Lie algebras in finite characteristics, the
irreducible representation corresponding to a given highest weight may
have a smaller dimension than the classical one
\cite{Lusztig,DobStAndrews,BKW,DobTruIrreg}. 

\medskip

Another feature can arise for $\uq{N}$ representations in the limit
when $q^l=1$: they can be {\em non integrable}, in the sense that 
the $\uq{N-1}\subset\uq{N}$ representations it contains may
become indecomposable.

\medskip

In the classical case and in the case of generic $q$, the
Gelfand--Zetlin basis for $\uq3$ irreducible representations
simultaneously  
diagonalizes the Cartan generators {\em and}
the quadratic Casimir operator $\csl2$ of $\uq2$ \cite{JimGZ}. 

When $q$ is a root of unity, some of the $\uq2$-modules involved in a
simple $\uq3$-module can be indecomposable with a non diagonalizable
action of the quadratic Casimir operator $\csl2$. 
If the definition of the \GZ
basis includes the requirement that $\csl2$ is diagonalized, then this
basis cannot exist for such a representation. If we consider the weaker
requirement that $\csl2$ is expressed 
in indecomposable blocks, then the \GZ
basis exists, as we will show. Indeed, in the limit when $q$ is a
root of unity of order $l$, the $\uq2$ representations 
of dimensions $l+d$ and $l-d$
have the same value of $\csl2$ and they are coupled in a single
indecomposable representation.

The signal that the \GZ basis without modification does not work for
non integrable 
irreducible restricted representations at roots of unity is given
by the fact that some denominators vanish in the coefficients. No
scale change can solve this problem. As explained in \cite{Abdatyp},
solution to cure the 
divergences is a suitable mixing
of states with the same quantum numbers. Since this mixing involves
zero or infinite coefficients, the correct way is to perform it at 
generic $q$ and to take the limit. The limit of {\em all}
matrix elements being zero or finite, we get a well-defined description of
the representation at $q^l=1$. 

With such a well-defined description, it is then possible to
exhibit the subrepresentation, in the cases when it exists. 

\medskip

The results of this paper may be summarized as follows: 

\begin{itemize}
\item In the classical case, or when $q$ is generic, the total $\uq2$
  representation corresponding to a given value of the Cartan element
  $h_1+2h_2$ that commutes with $\uq2$ is equivalent to the tensor
  product of two $\uq2$ irreducible representations, the rule being shown in
  Figure \ref{fig:sl2subrep}.  
  When $q^l=1$, this property remains true, but the tensor product now
  decomposes into indecomposable and irreducible representations. 
\item If we introduce the two transformations acting on
  Gelfand--Zetlin states (the definitions are given in the next section)
  \begin{eqnarray}
    S_1~: 
    && \Ket{p_{13}}{p_{23}}{p_{33}}{p_{12}}{p_{22}}{p_{11}}
    \longmapsto
    \Ket{p_{13}}{p_{23}}{p_{33}}{p_{22}+l}{p_{12}-l}{p_{11}} \;,
    \label{eq:S1}
  \end{eqnarray}
  \begin{eqnarray}
    S_2~: && V(p_{33}+l,p_{23},p_{13}-l) \longrightarrow
    V(p_{13},p_{23},p_{33}) \nonumber\\
    && \Ket{p_{33}+l}{p_{23}}{p_{13}-l}{p_{12}}{p_{22}}{p_{11}}
    \longmapsto
    \Ket{p_{13}}{p_{23}}{p_{33}}{p_{12}}{p_{22}}{p_{11}} \;,
    \label{eq:S2}
  \end{eqnarray}
  ($S_2$ is defined when $p_{33}+l>p_{23}>p_{13}-l$ only)
  then 
  \begin{itemize}
  \item If a state {\em and} its image by $S_1$ belong to the
    representation, then these two sates belong to the same 
    indecomposable $\uq2$ representation and should be redefined.
    The set of redefined states looks like a teepee or a tent,
    depending on the highest weight (See Figure \ref{fig:teepee}).
    This happens when the highest weight $\lambda$ is such that
    $\langle\lambda,\theta^\vee\rangle \ge l$, where $\theta$ is the
    longest root. 
  \item The image of $S_2$ is a subrepresentation. 
    This image is exactly the subrepresentation described in  
    \cite{DobTruIrreg}.
    This happens when the highest weight $\lambda$ is such that
    $\langle\lambda+\rho,\theta^\vee\rangle >l$     
    and when its image by the reflection with respect to the line
    $\langle\lambda+\rho,\theta^\vee\rangle =l$ 
    is also a dominant weight 
    ($\rho$ being the sum of fundamental weights). 
  \end{itemize}
\end{itemize}

The transformations $S_1$ and $S_2$ are particular cases of
transformations introduced in \cite{ACperiodicsln} 
for periodic representations, 
corresponding to {\em i)} symmetry among the \GZ indices of the same
line: permutations of these indices leave the coefficients invariant
{\em ii)} invariance under a translation by $l$ of a \GZ index. 
These symmetries become a problem for the restricted representations
we consider here. The mixing and normalization of states we 
introduce actually break them.
The transformation $S_2$ now defines an isomorphism from
$M_q(p_{33}+l,p_{23},p_{13}-l)$ to a subrepresentation of
$M_q(p_{13},p_{23},p_{33})$. 

\medskip

The structure of this paper is the following: in Section
\ref{sect:definitions}, we recall the definition of the quantum
enveloping algebra $\uq3$ and give the expression of the \GZ basis
for generic deformation parameter $q$. In Section
\ref{sect:regularization}, we propose a mixing of some states of the
\GZ basis that allows a well-defined limit when $q^l=1$. In Section
\ref{sect:reducibility}, the subrepresentation of the regularized
representation is exhibited, when it exists. Finally, some technical
expressions are given in Appendices \ref{app:indec}, \ref{app:e2f2}
for indecomposable $\uq2$ 
representations and for action within the set of redefined states. 
Final checks of finiteness of coefficients are made in Appendix
\ref{app:divergences}.

\section{\label{sect:definitions}Definitions}

Let $\uq3$ be the unital algebra generated by $e_i$, $f_i$ and $h_i$
($i=1,2$) with the relations
\begin{eqnarray}
  {} [h_i, e_j] & = & a_{ij} e_j 
 \qquad \qquad
 {} [h_i f_j]  = -a_{ij} f_j \;,\nonumber\\
 {} [e_i,f_j] & = & \delta_{ij} \frac{q^{h_i} - q^{-h_i}}{ q-q^{-1}}
 = \delta_{ij} [h_i] \;,\nonumber\\
 e_i^2 e_{i\pm 1} & - & (q+q^{-1}) e_i e_{i\pm 1} e_i + e_{i\pm
   1} e_i^2 = 0 \;,\nonumber\\
 f_i^2 f_{i\pm 1} & - & (q+q^{-1}) f_i f_{i\pm 1} f_i + f_{i\pm
   1} f_i^2 = 0 \;,
  \label{eq:sl3}
\end{eqnarray}
where 
$(a_{ij})=\left(
\begin{array}{cc} 2 & -1 \\ -1 & 2  \end{array}
\right)$ 
is the Cartan matrix of $sl(3)$. We define $q$-numbers by
$[x]\equiv\frac{q^x-q^{-x}}{q-q^{-1}}$. 

\medskip

Let $\alpha_1$, $\alpha_2$ be the simple roots, $\omega_1$,
$\omega_2$ the fundamental weights, $P=\ZZ\omega_1\oplus\ZZ\omega_2$
the weight lattice and $\alpha_1^{\vee}$,
$\alpha_2^{\vee}$ the coroots, with 
$\langle \omega_i,\alpha_j^{\vee}\rangle = \delta_{ij}$.
The longest root is $\theta=\alpha_1+\alpha_2$. The sum of the
fundamental weights (equal to half the sum of positive roots) is
$\rho=\alpha_1+\alpha_2$. 

\medskip

We will later be interested in the root of unity case. {\em Let $l>2$ be an
odd integer.} When $q^l=1$, $\left(q^{h_i}\right)^{l}$ is central and
we will add the relations corresponding to the restricted specialization
\begin{eqnarray}
  &&  e_1^l=e_2^l=e_3^l=0 \;,\qquad \mbox{with} 
  \qquad e_3=e_1e_2-q^{-1}e_2e_1
  \;,\nonumber\\
  &&  f_1^l=f_2^l=f_3^l=0 \;,\qquad \mbox{with} 
  \qquad f_3=f_2f_1-qf_1f_2
  \;,\nonumber\\
  &&  q^{2h_i}_i = 1  \;.
  \label{eq:ql=1}
\end{eqnarray}
These relations define a co-ideal with respect to the Hopf structure,
so that quotienting by them leads to a Hopf algebra. We do not
introduce the divided powers of the generators. The generators $e_i$,
$f_i$ and $k_i$ and the relations (\ref{eq:sl3},\ref{eq:ql=1})
actually define a finite dimensional Hopf subalgebra of the usual
restricted specialization. As proved in \cite{Lusztig}, the study of finite
dimensional representations of the restricted specialization can be
reduced to the study of those of the finite subalgebra.

The finite dimensional irreducible representations are labeled by
integral dominant weights
$\lambda=\lambda_1\omega_1+\lambda_2\omega_2$, with
$\lambda_i\in\ZZ_+$. 
We can limit the study to $0\le \lambda_i < l$
since translations of the highest weight by multiples of $l\omega_i$ 
provide equivalent representations (strictly speaking, the
representations are only equivalent as representations of the algebra
generated by $e_i$, $f_i$ and $q^{h_i}$, $(i=1,2)$; a global
translation of the weights is however the only difference).

\medskip

As in the case of affine Lie algebras, a representation $M$ is called
integrable (see, e.g. \cite{LasLecThi}) if 
\begin{enumerate}
\item $M=\bigoplus_{\Lambda\in P} M_\Lambda$, i.e. $M$ is the direct
  sum of its weight spaces (common eigenspaces of the Cartan
  generators), the weights being integral (belonging to the weight
  lattice $P$), 
\item $\dim M_\Lambda < \infty$, i.e. each weight space has a finite
  dimension, 
\item $M$ decomposes into a direct sum of finite dimensional
  representations of the $\cU_q(sl(2))$ subalgebras generated by
  $e_i$, $f_i$, $q^{h_i}$, $q^{-h_i}$ for $i=1,2$.
\end{enumerate}

In the case we consider, all irreducible representations have a finite
dimension, since, at roots of unity, the quantum algebra is a finite
dimensional module over its centre. The first two requirements for
integrability are hence always satisfied for irreducible
representations. As we will see, the third one is not always satisfied
since $M$ may contain indecomposable representations of its quantum
subalgebras. 

\medskip

For generic $q$, any finite dimensional irreducible representation
can be described using the Gelfand--Zetlin basis. Let
$M_q(p_{13},p_{23},p_{33})$ be the representation
with highest weight 
$(\lambda_1,\lambda_2)=(p_{13}-p_{23}-1,p_{23}-p_{33}-1)$, 
(the eigenvalues of $h_1$ and $h_2$ on the highest weight vector). 
It acts on the vector space $V(p_{13},p_{23},p_{33})$ of dimension 
\begin{equation}
d(p_{13},p_{23},p_{33}) =\frac12 (p_{13}-p_{23})(p_{23}-p_{33})(p_{13}-p_{33})
  \label{eq:dim}
\end{equation}
and spanned by vectors 
\begin{equation}
  \Ket{p_{13}}{p_{23}}{p_{33}}{p_{12}}{p_{22}}{p_{11}}
  \label{eq:ket}
\end{equation}
with $p_{ij}\in\ZZ$, such that
\begin{eqnarray}
  p_{13} &\ge & p_{12} > p_{23} \ge p_{22} > p_{33} \;,\nonumber\\
  && p_{12} \ge p_{11} > p_{22} \;.
  \label{eq:triangineq}
\end{eqnarray}
All the $p_{ij}$ are defined up to an overall constant. Only
differences are involved in the matrix elements. We use
$p_{ij}=h_{ij}-i$ instead of the standard $h_{ij}$ to make more
explicit the symmetries among the indices of the same line. 
The first line of indices is constant for a given representation. We
will sometimes omit it, when no confusion is possible and when it is
the same as in (\ref{eq:ket}). 

\medskip

The representations $M_q(p_{13},p_{23},p_{33})$ described here are
actually in one-to-one correspondence with the classical ones. To get
all the $2^{\expo{rank}}=4$ inequivalent 
representations corresponding to a classical, one can add to the
index $p_{11}$, or to the index $p_{12}$, or to both of
them, the constant $i\pi/\ln q$ (or $l/2$ if $q^l=1$).

\medskip

At generic $q$ as well as in the classical case, the \GZ basis
expresses the $\uq3$ representation as a direct sum of $\uq2$
irreducible representations (corresponding to fixed values of $p_{12}$
and $p_{22}$). By $\uq2$, we will always mean the
subalgebra of $\uq3$ generated by $e_1$, $f_1$, $k_1$. 

\medskip

As in \cite{DobTruIrreg}, we shall depict the set of \GZ state in a
three dimensional pyramid, with one point for each vector of the
basis. The horizontal coordinates $x,y$ 
are simply the values of the orthogonal Cartan
elements $h_1$ and $h_1+2h_2$. The third coordinate $z$ 
starts form $0$ and
increases when the dimension $p_{12}-p_{22}$ 
of the $\uq2$ representation decreases.
\begin{eqnarray}
  x &=& 2p_{11} - (p_{12}+p_{22}) -1 \;,\nonumber \\
  y &=& 3(p_{12}+p_{22}) -2(p_{13}+p_{23}+p_{33})-1 \;,\nonumber\\
  z &=& \min(p_{13}-p_{12},p_{23}-p_{33}-1) \;.
  \label{eq:coord}
\end{eqnarray}

The actions of the generators on the \GZ basis are given by
\begin{eqnarray}
  h_1 |p\rangle &=& {\big(2p_{11}-(p_{12}+p_{22})-1\big)} |p\rangle
  \;,\nonumber \\
  h_2 |p\rangle &=&
  {\big(2(p_{12}+p_{22})-p_{11}-(p_{13}+p_{23}+p_{33})-1\big)} 
  |p\rangle \;,
  \label{eq:GZ0}
\end{eqnarray}
\begin{eqnarray}
  f_1 |p\rangle &=& \Big([p_{12}-p_{11}+1][p_{11}-p_{22}-1]\Big)^{1/2}
  \ket{p_{12}}{p_{22}}{p_{11}-1}  
  \;,\nonumber \\
  e_1 |p\rangle &=& \Big([p_{12}-p_{11}][p_{11}-p_{22}]\Big)^{1/2}
  \ket{p_{12}}{p_{22}}{p_{11}+1}  
  \;,
  \label{eq:GZ1}
\end{eqnarray}
\begin{eqnarray}
  f_2 |p\rangle &=& \left(\frac {P_1 P_2}{P_3}(1,2;p)
  \right)^{1/2}|p_{12}-1\rangle 
  + \left(\frac {P_1 P_2}{P_3}(2,2;p)
  \right)^{1/2}|p_{22}-1\rangle 
  \;,\nonumber\\
  e_2 |p\rangle &=& \left(\frac 
  {P_1 P_2}{P_3}(1,2;p_{12}+1)
  \right)^{1/2}|p_{12}+1\rangle 
  + \left(\frac {P_1 P_2}{P_3}(2,2;p_{22}+1)
  \right)^{1/2}|p_{22}+1\rangle 
  \;,
  \label{eq:GZ2}
\end{eqnarray}
where 
\begin{eqnarray}
  P_1(1,2;p) &=& [p_{13}-p_{12}+1][p_{12}-p_{23}-1][p_{12}-p_{33}-1]
  \;,\nonumber\\
  P_1(2,2;p) &=& [p_{13}-p_{22}+1][p_{23}-p_{22}+1][p_{22}-p_{33}-1]
  \;,\nonumber\\
  P_2(1,2;p) &=& [p_{12}-p_{11}] \;,\nonumber\\
  P_2(2,2;p) &=& [p_{11}-p_{22}] \;,\nonumber\\
  P_3(1,2;p) &=& [p_{12}-p_{22}][p_{12}-p_{22}-1] \;,\nonumber\\
  P_3(2,2;p) &=& [p_{12}-p_{22}][p_{12}-p_{22}+1] \;,
  \label{eq:P1P2P3}
\end{eqnarray}
where $p$ stands for the set of indices $p_{ij}$, and where 
$p_{ij}\pm 1$ in an argument shows the modified index only. The two
first arguments $i,j$ of the coefficients $P_{\alpha}$ indicate which
$p_{ij}$ is changed. 

\medskip

For generic $q$, the $q$-integers involved in the coefficients vanish
only at zero argument. Vanishing denominators are compensated by two
vanishing numerators. 

\medskip

When $q$ goes to a primitive $l$-th root of one,
with $l$ odd, 
the $q$-integer $[n]$ goes to zero iff $n$ is a multiple of $l$. 
For this reason, 
new zeroes arise in the denominator when $p_{13}-p_{33}-2\ge l$, i.e.\  when
the highest weight satisfies 
$\langle\lambda,\theta^\vee\rangle = \lambda_1+\lambda_2\ge l$. 
These new zeroes are
generally not compensated in the numerator. The previously defined \GZ
basis is then not well-defined in this case.

When 
$\langle\lambda,\theta^\vee\rangle = p_{13}-p_{33}-2 < l$, 
the representation is correctly described
by the \GZ basis. The $\uq2$ representations it involves are
completely reducible into irreducible representations of 
dimension less than $l$.

The remaining case is then $p_{13}-p_{33} -2\ge l$ and still 
$p_{13}-p_{23}\le l$ and $p_{23}-p_{33}\le l$.  It is the aim of
Section \ref{sect:regularization} to get a well-defined description of
this case.

\section{\label{sect:regularization}Regularization}

We still consider a representation $M_q(p_{13},p_{23},p_{33})$ of 
$\uq3$ at {\em generic} $q$.
Let $l>2$ be an odd integer. When  $p_{13}-p_{33}>l$
(and $p_{13}-p_{23}\le l$, $p_{23}-p_{33}\le l$), 
in prevision of the case $q^l=1$, 
we perform the following transformation that depends on $l$.

Let us consider $\ket{p_{12}}{p_{22}}{p_{11}}$ with
$p_{12}-p_{22}>l$. 
When both (\ref{eq:triangineq}) and 
\begin{eqnarray}
  p_{13} &\ge & p_{22}+l > p_{23} \ge p_{12}-l > p_{33} \;,\nonumber\\
  && p_{22}+l \ge p_{11} > p_{12}-l
  \label{eq:triangineq2}
\end{eqnarray}
are satisfied, i.e.\  if the image of $\ket{p_{12}}{p_{22}}{p_{11}}$ by
$S_1$ defined by (\ref{eq:S1}) also belongs to the representation, 
we define
\begin{equation}
  \left( 
    \begin{array}{c}
      \ket{p_{12}}{p_{22}}{p_{11}}' \\
      \ket{p_{22}+l}{p_{12}-l}{p_{11}}' \\
    \end{array}
  \right) 
  =
  \left(
    \begin{array}{cc}
      [l]^{1/2} & 0 \\
      \left(\fracdis{[l-1]}{[l+1][l]}\right)^{1/2} & \fracdis1{[l]^{1/2}} 
    \end{array}
  \right)
  \left( 
    \begin{array}{c}
      \ket{p_{12}}{p_{22}}{p_{11}} \\
      \ket{p_{22}+l}{p_{12}-l}{p_{11}} \\
    \end{array}
  \right) \;.
  \label{eq:transf}
\end{equation}
This transformation is inspired by that introduced in \cite{Abdatyp}.
As in this paper, the transformation matrix has determinant $1$ and its 
eigenvalues go to $0$ and $\infty$ in the limit when $q^l=1$. 
In the following, we keep the primed states and forget the
corresponding unprimed states. 

The set of \GZ states such that their image by $S_1$ is still a \GZ
state (i.e.\  satisfying both (\ref{eq:triangineq}) and 
(\ref{eq:triangineq2}))
is displayed (on the hexagonal basis of the pyramid) 
in Figure \ref{fig:teepee}. 
It looks like a teepee or a tent, depending on the values of $p_{i3}$. 
It includes 
both the redefined states and the 
$\uq2$ representations 
invariant under $S_1$, with dimension $p_{12}-p_{22}=l$, that are not
redefined.


\begin{figure}
  \setlength{\unitlength}{0.015cm}
  \begin{picture}(655,420)(-200,0)
    \drawline(500,35)(325,370)(240,10)
    (155,160)(315,388)
    \drawline(155,160)(0,115)(55,20)
    (405,0)(555,55)(455,165)(430,170)
    \drawline(315,388)(325,370)(375,155)
    (365,55)(240,10)
    \drawline(375,155)(240,10)
    \drawline(375,155)(500,35)(365,55)
  \end{picture}
  \caption{Teepee\label{fig:teepee}}
\end{figure}

\subsection{Indecomposable $\uq2$ subrepresentations}

\begin{figure}[htbp]
  \setlength{\unitlength}{0.01in}%
  \begin{center}
    \leavevmode
    \begin{picture}(600,225)(80,400)
      \thicklines
      \multiput(300,550)(50,0){4}{\circle*{7}}
      \multiput(150,500)(50,0){10}{\circle*{7}}
      \put(300,550){\line( 1, 0){150}}
      \put(150,500){\line( 1, 0){450}}
      \multiput(300,550)(50,0){4}{\vector(1, -1){48}}
      \multiput(300,550)(50,0){4}{\vector(-1, -1){48}}
      \multiput(150,500)(50,0){6}{\vector(1, 0){48}}
      \multiput(500,500)(50,0){2}{\vector(1, 0){48}}
      \multiput(300,550)(50,0){3}{\vector(1, 0){48}}
      \multiput(350,550)(50,0){3}{\vector(-1, 0){48}}
      \multiput(350,500)(50,0){6}{\vector(-1, 0){48}}
      \multiput(200,500)(50,0){2}{\vector(-1, 0){48}}
      \put(450,500){\vector(2, -1){22}}
      \put(480,480){\makebox(0,0)[lb]{{$\scriptstyle 0$}}}
      \put(300,500){\vector(-2, -1){22}}
      \put(270,480){\makebox(0,0)[lb]{{$\scriptstyle 0$}}}
      \put(353,450){\vector(1, 0){41}}
      \put(353,450){\vector(1, -1){41}}
      \put(347,450){\vector(-1, 0){41}}
      \put(347,450){\vector(-1, -1){41}}
      \put(365,410){\makebox(0,0)[lb]{\raisebox{0pt}[0pt][0pt]{$e_1$}}}
      \put(322,410){\makebox(0,0)[lb]{\raisebox{0pt}[0pt][0pt]{$f_1$}}}
    \end{picture}
  \end{center}
  \caption{Indecomposable $\cU_q(sl(2))$ representation with $q^7=1$ and
    $p_{12}-p_{22}=10$.} 
  \label{fig:indec}
\end{figure}
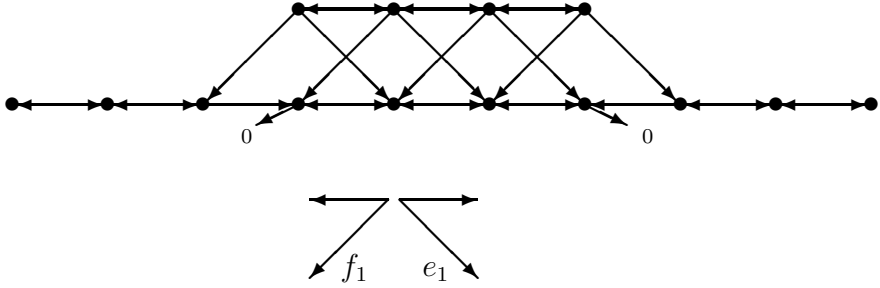

The two $\uq2$ representations corresponding to $p_{12}\ge p_{11}>p_{22}$
(of dimension $p_{12}-p_{22}$) and $p_{22}+l\ge p_{11}>p_{12}-l$ (of
dimension $2l-(p_{12}-p_{22})$) are gathered into a sum of dimension
$2l$. 
In the limit when $q^l=1$, this sum becomes indecomposable. It
is described in Figure \ref{fig:indec}. 
The actions of the generators $e_1$ and $f_1$ 
on this indecomposable
representation are given in Appendix \ref{app:indec}. Initially, $e_1$
and $f_1$ induced only moves along the $x$ direction. Now, the primed
states replace the unprimed states inside the teepee, and $e_1$
and $f_1$ induce moves along the $x$ direction {\em and} possibly
shortcut down in the $z$ direction. In the extreme case, this shortcut can
lead directly from top to bottom. 

The quadratic Casimir operator of $\uq2$ acts on the space spanned by 
$\ket{p_{12}}{p_{22}}{p_{11}}'$ and 
$\ket{p_{22}+l}{p_{12}-l}{p_{11}}'$ as the non diagonalizable matrix
\begin{eqnarray}
  \csl2 
  =
  \left(
    \begin{array}{cc}
      q^{p_{12}-p_{22}}+q^{p_{22}-p_{12}} & 0 \\
      i (q-q^{-1})^2 
      [p_{12}-p_{22}]
      &
      q^{p_{12}-p_{22}}+q^{p_{22}-p_{12}} 
    \end{array}
  \right) \;.
  \label{eq:c2}
\end{eqnarray}

\begin{figure}
\newcounter{ve}
\setlength{\unitlength}{0.01in}%
\begin{picture}(210,300)(120,390)
  \put(440,640){\makebox(0,0)[lb]{ $j_1$ }}
  \put(400,630){\vector(1,0){100}}
  \put(500,630){\vector(-1,0){100}}
  \multiput(400,630)(0,-3){10}{\circle*{0.2}}
  \multiput(500,630)(0,-3){10}{\circle*{0.2}}
  \put(370,600){\vector(-2,-3){75}}
  \put(295,487.5){\vector(2,3){75}}
  \multiput(370,600)(3,0){10}{\circle*{0.2}}
  \multiput(295,487.5)(3,0){10}{\circle*{0.2}}
  \put(300,540){\makebox(0,0)[lb]{ $j_2$ }}
  \put(325,480){\vector(1,0){250}}
  \put(575,480){\vector(-1,0){250}}
  \put(425,465){\makebox(0,0)[lb]{ $j_1\otimes j_2$ }}
  \setcounter{ve}{100}
  \multiput(400,600)(-25,-37.5){5}{\line(1, 0){
      \arabic{ve}}\addtocounter{ve}{50}}
  \multiput(300,412.5)(10,0){30}{\line(1,0){5}}
  \multiput(275,412.5)(-20,-30){2}{\line(-2,-3){10}}
  \multiput(625,412.5)(20,-30){2}{\line(2,-3){10}}
  \setcounter{ve}{3}
  \multiput(400,600)(-25,-37.5){5}{
    \multiput(0,0)(50,0){\arabic{ve}}{\circle*{4}}
    \addtocounter{ve}{1}}
  \setcounter{ve}{2}
  \multiput(425,562.5)(-25,-37.5){4}{
    \multiput(0,0)(50,0){\arabic{ve}}{\circle{8}}
    \addtocounter{ve}{1}}
  \setcounter{ve}{1}
  \multiput(450,525)(-25,-37.5){3}{
    \multiput(0,0)(50,0){\arabic{ve}}{\circle{12}}
    \addtocounter{ve}{1}}
  \put(400,600){\line( -2, -3){120}}
  \put(500,600){\line( 2, -3){120}}
\end{picture}
\caption{Rule for $\cU_q(sl(2))$ subrepresentation.\label{fig:sl2subrep}}
\end{figure}
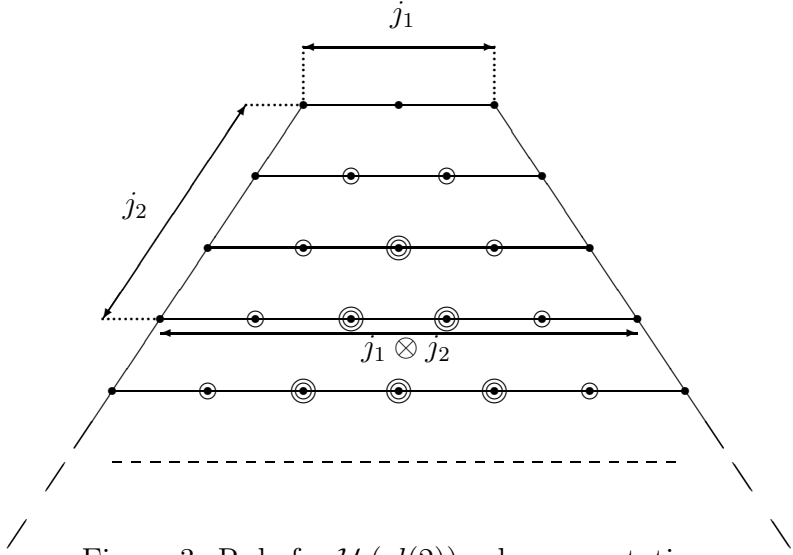

To summarize, the $\uq2$ modules with the same value of $h_1+2h_2$
that would have the same value of the quadratic Casimir when $q^l=1$
are pairwise 
coupled in a single indecomposable representation. The same thing
happens in the fusion rule of restricted irreducible representations
of $\uq2$ \cite{PS,Kel,Suter}. The total $\uq2$ representation
corresponding 
to a given value of $h_1+2h_2$ is actually equivalent, as in the
classical or generic case, to the tensor product of two irreducible
representations. If the value of
$h_1+2h_2$ is higher or equal to $p_{13}-2p_{23}+p_{33}+2$
(corresponding to the classical $sl(2)$ representation with the highest
dimension $p_{13}-p_{33}-1$), the total $\uq2$ representation
corresponding to this value is equivalent to 
\begin{equation}
  j_1 \otimes j_2 
  \qmbox{with}
  \cases{
    j_1 = \frac12 (p_{13}-p_{23}-1) \cr
    j_2 = \frac16 (p_{13}+p_{23}-2p_{33}-1-(h_1+2h_2))}
  \label{eq:tensor1}
\end{equation}
as shown on Figure \ref{fig:sl2subrep} which is easier to understand
than the formula.
In the case where the value of
$h_1+2h_2$ is lower or equal to $p_{13}-2p_{23}+p_{33}+2$, the total
$\uq2$ representation is equivalent to
\begin{equation}
  j_1 \otimes j_2 
  \qmbox{with}
  \cases{
  j_1 = \frac12 (p_{23}-p_{33}-1) \cr
  j_2 = \frac16 (h_1+2h_2 - (-2p_{13}+p_{23}+p_{33}-1))}
\label{eq:tensor2}
\end{equation}
i.e.\  the same as before, but starting from the opposite edge of the
hexagon. 

\subsection{Other effects of the regularization}

In the limit when $q^l=1$, the coefficients that would involve
fractions like $\frac{\displaystyle[0]}{\displaystyle[l]}$ remain
$0$. Although $[l]$ is zero at 
$q^l=1$, there is no ambiguity in such limits. This means in
particular that the states that are classically forbidden (those that
would not respect the triangular inequalities (\ref{eq:triangineq})) 
remain forbidden at $q^l=1$. The vector space
$V(p_{13},p_{23},p_{33})$ on which the representation acts at $q^l=1$
is the same as in the classical case or at generic $q$. (As we will
see in Section \ref{sect:reducibility}, the so-obtained representation
is however sometimes not irreducible, and we will be led to take
quotients.)  

Since the redefinition (\ref{eq:transf}) contain coefficients that
diverge in the limit $q^l=1$, we have to check carefully the behaviour
of the regularized representation on redefined states and when crossing 
the boundary of the teepee, the
domain that contains the redefined states. 

\begin{itemize}
\item In Appendix
  \ref{app:indec}, the actions of $f_1$ and $e_1$ are explicitly
  given. The coefficients are finite. As explained before, they
  describe well-defined indecomposable representations of $\uq2$. 
\item In Appendix \ref{app:e2f2}, the actions of $e_2$ and $f_2$ within
  the teepee are computed. The coefficients are also finite. The
  diverging 
  coefficient in (\ref{eq:transf}) essentially enhances differences
  of coefficients that have the same limits.
\item In Appendix \ref{app:divergences}, a compendium of all possible
  sources of divergences is presented. In each case, the reason why
  the divergence disappears is briefly explained.

  The locations where the generators can 
  lead from a non redefined state to a redefined one, or vice-versa,  
  is the set of \GZ states satisfying one of the following
  equations:  
  \begin{itemize}
  \item Left roof: $p_{22}=p_{13}-l$,
  \item Right roof: $p_{12}=p_{33}+l+1$,
  \item Front ``entrance'': $p_{22}=p_{11}-l$,
  \item Back ``entrance'': $p_{12}=p_{11}+l-1$,
  \item $l$-dimensional $\uq2$ modules: $p_{12}-p_{22}=l$.
  \end{itemize}
  The adjectives ``left'', ``right'', ... refer to Figure
  \ref{fig:teepee}. 
  The first four cases are the boundaries of the teepee within the
  pyramid. The last case corresponds to non redefined $l$-dimensional
  modules. The boundaries are defined as belonging to the teepee. 
  Note that the front and back roofs are boundaries of the pyramid,
  not boundaries of the teepee in the pyramid.
\end{itemize}
  
After the regularization defined by (\ref{eq:transf}), all the
coefficients then remain finite or go to zero the limit where $q^l=1$. 
A representation $M^{\expo{reg}}_q(p_{13},p_{23},p_{33})$ at $q^l=1$ is then
obtained by defining the action of 
the generators using the limit of these coefficients. These
coefficients being finite, they indeed define elements of
$\mbox{End}(V(p_{13},p_{23},p_{33}))$.
Moreover, these elements satisfy the commutation relations of $\uq3$
at $q^l=1$, since these relations are continuous functions of the
coefficients.

\section{\label{sect:reducibility}Reducibility}

The regularized representation 
$M^{\expo{reg}}_q(p_{13},p_{23},p_{33})$  at $q^l=1$ 
is not always irreducible. We recall that we consider
$p_{13}-p_{33}>l$ (Otherwise, nothing new happens with respect to the
generic case).

If $p_{23}$ is equal to $p_{13}-l$ or to $p_{33}+l$,
$M^{\expo{reg}}_q(p_{13},p_{23},p_{33})$  is irreducible.
Otherwise, i.e.\  if $\min(p_{13}-p_{23},p_{23}-p_{33})<l$,
the application $S_2$ defined in equation (\ref{eq:S2}) from the
vector space $V(p_{33}+l,p_{23},p_{13}-l)$ to the vector space
$V(p_{13},p_{23},p_{33})$ is a morphism from the representation 
$M_q(p_{33}+l,p_{23},p_{13}-l)$ to the representation
$M^{\expo{reg}}_q(p_{13},p_{23},p_{33})$. Its image is isomorphic to 
$M_q(p_{33}+l,p_{23},p_{13}-l)$, and is a subrepresentation of 
$M^{\expo{reg}}_q(p_{13},p_{23},p_{33})$. It is easy to check that 
\begin{itemize}
\item None of the redefined state belongs to this image.
\item No action of the generators connects directly this image to the
  set of redefined states.
\end{itemize}
That this image really decouples can then be seen using only equations
(\ref{eq:GZ1},\ref{eq:GZ2}). The factors $[p_{12}-p_{33}\,(-1)]^{1/2}$ and 
$[p_{13}-p_{22}\,(+1)]^{1/2}$, that vanish for $p_{12}=p_{33}+l\;(+1)$ and
$p_{22}=p_{13}-l\;(-1)$, respectively, are enough.

The representation $M^{\expo{reg}}_q(p_{13},p_{23},p_{33})$ is then
the direct sum of 
the two subrepresentations respectively characterized by 
$\max(p_{12}-p_{33},p_{13}-p_{22}+1)\ge l$ and 
$\max(p_{12}-p_{33},p_{13}-p_{22}+1)< l$. 

The first one, equivalent
to $M_q(p_{33}+l,p_{23},p_{13}-l)$, has then a classical counterpart. 
In \cite{DobTruIrreg}, this subrepresentation 
is identified with  the 
$\min(p_{33}-p_{23}+l,p_{23}-p_{13}+l)$
top layers of the pyramid. 

The second one, that contains
all the redefined states, has no classical analogue. It corresponds in
\cite{DobTruIrreg} to the 
$p_{13}-p_{33}-l$ 
bottom layers of the pyramid (Its height is the same as that of the
teepee). 
We denote it by $M^{\expo{quot}}_q(p_{13},p_{23},p_{33})$ as it is the
quotient of $M^{\expo{reg}}_q(p_{13},p_{23},p_{33})$ by 
$S_2(M_q(p_{33}+l,p_{23},p_{13}-l))$.
Its dimension is
$d(p_{13},p_{23},p_{33})-d(p_{33}+l,p_{23},p_{13}-l)$. 

These two summands are themselves irreducible. The reducibility of one
of the summands would require more singular vectors in the Verma
module with the same highest weight as $M_q(p_{13},p_{23},p_{33})$
than found in \cite{DobStAndrews}.

\section{An interesting case: flat representations}

An interesting case is provided by the flat representations,
i.e.\  those for which the weight multiplicities are at most $1$. They
correspond to parameters such that $p_{13}-p_{33}=l+1$. 

In this case,
no state needs being redefined, since the teepee reduces to one single
line with $l$ points 
(with $p_{12}=p_{13}$ and $p_{22}=p_{33}+1$). 
These representations are actually integrable in the sense given in
Section \ref{sect:definitions}, since
$\langle\lambda,\theta^\vee\rangle = l-1$, the maximum value for
integrable representations. 

If $p_{23}$ is equal
to $p_{13}-1=p_{33}+l$ or to $p_{13}-l=p_{33}+1$, then
$M_q(p_{13},p_{23},p_{33})$ itself is flat and irreducible.
Otherwise, the flat irreducible representation is, as explained before, 
$M_q(p_{13},p_{23},p_{33})/S_2(M_q(p_{13}-1,p_{23},p_{33}+1))$, of
dimension $d(p_{13},p_{23},p_{33})-d(p_{13}-1,p_{23},p_{33}+1)$.

The states of this quotients are then the \GZ states satisfying the
usual triangular inequalities (\ref{eq:triangineq}) and 
\begin{equation}
  p_{12}=p_{13} \qmbox{or} p_{22}=p_{33}+1=p_{13}-l \;.
  \label{eq:flat1}
\end{equation}

The existence and dimension of these representations were known from the
character formulas \cite{Lusztig,DobStAndrews}.

\medskip

The flat irreducible representations of $\uq3$ were described in
\cite{Asutroisq} as  
quotients of singular limits of flat periodic representations of
dimension $l^2$. They
were also obtained in \cite{ACperiodicsln,AACuqsln} 
within the \GZ basis, but
with a different prescription that we recall now .  
Consider the vector space spanned by the vectors
\begin{equation}
  \Ket{\bar p_{13}}{\bar p_{23}}{\bar p_{23}-1}
  {\bar p_{12}}{\bar p_{23}}{\bar p_{11}}
  \label{eq:ketflat}
\end{equation}
with $\bar p_{23}+2l\ge \bar p_{13} > \bar p_{23}+l$ and where  
$\bar p_{22}=\bar p_{23}$ is frozen. 


With the triangular inequalities (\ref{eq:triangineq}), this defines a
$\uq3$ representation with a triangular set of weights (hence
flat) of dimension $d_0=d(\bar p_{13},\bar p_{23},\bar
p_{23}-1)$. This representation is not irreducible and splits into
four subrepresentations obtained as follows:
\begin{center}
  \begin{tabular}{p{4cm}p{11.5cm}}
    $\bar p_{13} \ge \bar p_{12} > \bar p_{23}+l$ 
    and 
    $\bar p_{12} \ge \bar p_{11} > \bar p_{23}+l$ 
    & Flat triangular representation (left corner) of dimension
    $d_1=d(\bar p_{13},\bar p_{23}+l,\bar p_{23}+l-1)
    =\frac12(\bar p_{13}-\bar p_{23}-l)(\bar p_{13}-\bar p_{23}-l+1)$. 
    \\[4mm]
    $\bar p_{13} \ge \bar p_{12} > \bar p_{23}+l$ 
    and 
    $\bar p_{12}-l \ge \bar p_{11} > \bar p_{23}$ 
    & Another flat triangular representation 
    (right corner), of dimension 
    $d_2=d(\bar p_{13}-l,\bar p_{23}-l,\bar p_{23}-l-1) = d_1$.  
    \\[4mm]
    $\bar p_{13} -l \ge \bar p_{12} > \bar p_{23}$ 
    and 
    $\bar p_{12} \ge \bar p_{11} > \bar p_{23}$ 
    & Another flat triangular representation 
    (bottom corner), of dimension 
    $d_3=d(\bar p_{13}-l,\bar p_{23}-l,\bar p_{23}-l-1) = d_1$.  
    \\[4mm]
    $\bar p_{13} \ge \bar p_{12} > \bar p_{23}$ 
    and 
    $\bar p_{12} \ge \bar p_{11} > \bar p_{12}-l$ 
    and 
    $\bar p_{23}+l \ge \bar p_{11} > \bar p_{23}$ 
    & The flat
    hexagonal representation of dimension $d_0-d_1-d_2-d_3=d_0-3d_1$.
  \end{tabular}
\end{center}

This description is linked to the \GZ formalism of this paper by 
the identification $\bar p_{13}=p_{23}+l$, $\bar p_{23}=p_{13}-l$,
$\bar p_{33}=p_{33}=p_{13}-l-1$. A transformation inspired by both
$S_1$ and $S_2$ relates them, namely
\begin{eqnarray}
  \Ket{p_{13}}{p_{23}}{p_{33}}{p_{12}}{p_{22}}{p_{11}}
  &\longmapsto&
  \Ket{p_{23}+l}{p_{13}-l}{p_{33}}
  {p_{12}}{p_{22}}{p_{11}} 
  \nonumber\\[2mm]
  && \mbox{if $p_{22}=p_{13}-l$,}  \nonumber\\
  \Ket{p_{13}}{p_{23}}{p_{33}}{p_{12}}{p_{22}}{p_{11}}
  &\longmapsto&
  \Ket{p_{23}+l}{p_{13}-l}{p_{33}}
  {p_{22}+l}{p_{12}-l}{p_{11}} 
  \nonumber\\[2mm]
  && \mbox{if $p_{12}=p_{13}$.}  
  \label{eq:newS}
\end{eqnarray}

\section{\label{sect:conclusion}Conclusion}

The description of $\uq{N}$ representations with the Gelfand--Zetlin
basis will necessitate the knowledge, within the \GZ basis, of some
indecomposable $\uq{N-1}$ representations, probably those involved in
the decomposition of tensor products of irreducible ones \cite{FGP}. 
The representations will be
built by collecting the $\uq{N-1}$ representations that would have the
same values of the Casimir operators $\csl{N-1}^{(i)}$ in the limit
$q^l=1$. The use of transformations that generalize $S_1$ will help in
characterizing them. This being done, analogues of $S_2$ will provide
the subrepresentations.

The restricted representations we have described here can be used as
explained in \cite{AACuqsln} to build special kinds of partially
periodic (unrestricted) 
irreducible representations of $\uq{N}$ with $N>3$. 

One could also wonder whether the periodic indecomposable representations
of $\uq{2}$ of dimension $2l$ that arise in the fusion of periodic
irreducible representations \cite{Akinetic} may also 
appear in some $\uq{3}$ irreducible representations. 

\bigskip
{\bf Acknowledgements:} the author wishes to thank B. Abdesselam and
  F. Barbarin for discussions and communication of preprints.

\appendix

\section{Indecomposable $\uq2$ subrepresentations\label{app:indec}}

An indecomposable representation of dimension $2l$ is made from 
the two $\uq2$ representations corresponding to $p_{12}\ge p_{11}>p_{22}$
(of dimension $p_{12}-p_{22}$) and $p_{22}+l\ge p_{11}>p_{12}-l$ (of
dimension $2l-(p_{12}-p_{22})$). The states that have a common
$p_{11}$ are mixed as explained in (\ref{eq:transf}) and the limit
$q^l=1$ is taken.
The generators $e_1$ and $f_1$ act on it as 
\begin{eqnarray}
  && f_1 \ket{p_{12}}{p_{22}}{p_{11}} 
  = \Big([p_{12}-p_{11}+1][p_{11}-p_{22}-1]\Big)^{1/2}
  \ket{p_{12}}{p_{22}}{p_{11}-1}  
  \nonumber\\
  &&  \qquad\qmbox{for} p_{11}>p_{22}+l+1     \qmbox{or}
    p_{11}\le p_{12}-l
  \;,\nonumber \\
  && f_1 \ket{p_{12}}{p_{22}}{p_{11}}' 
  = \Big([p_{12}-p_{11}+1][p_{11}-p_{22}-1]\Big)^{1/2}
  \ket{p_{12}}{p_{22}}{p_{11}-1}'  
  \nonumber\\
  &&  \qquad\qmbox{for}  p_{22}+l\ge p_{11}> p_{12}-l+1
  \;,\nonumber \\
  && f_1 \ket{p_{12}}{p_{22}}{p_{22}+l+1} =
  \Big([p_{12}-p_{22}-l]\Big)^{1/2}
  \ket{p_{12}}{p_{22}}{p_{22}+l}' \qquad \mbox{(instead of $0$)}
  \;,\nonumber\\
  && f_1 \ket{p_{12}}{p_{22}}{p_{12}-l+1}' = 0
  \;,\nonumber\\
  && f_1 \ket{p_{22}+l}{p_{12}-l}{p_{12}-l+1}' = 
  \Big(-[p_{12}-p_{22}-l]\Big)^{1/2}
  \ket{p_{12}}{p_{22}}{p_{12}-l} 
  \nonumber\\
  &&\qquad\qquad \mbox{(instead of $0$)}
  \;,\nonumber\\
  && f_1 \ket{p_{22}+l}{p_{12}-l}{p_{11}}' = 
  \Big([p_{12}-p_{11}+1][p_{11}-p_{22}-1]\Big)^{1/2}
  \ket{p_{22}+l}{p_{12}-l}{p_{11}-1}' 
  \nonumber\\
  &&\qquad \qquad
  + \Big(-[p_{12}-p_{11}+1][p_{11}-p_{22}-1]\Big)^{-1/2} [p_{12}-p_{22}]
  \ket{p_{12}}{p_{22}}{p_{11}-1}  
  \nonumber\\
  &&  \qquad\qmbox{for}  p_{22}+l\ge p_{11}> p_{12}-l+1 \;.
  \label{eq:indec-f1}
\end{eqnarray}

\begin{eqnarray}
  && e_1 \ket{p_{12}}{p_{22}}{p_{11}} 
  = \Big([p_{12}-p_{11}][p_{11}-p_{22}]\Big)^{1/2}
  \ket{p_{12}}{p_{22}}{p_{11}+1}  
  \nonumber\\
  &&  \qquad\qmbox{for}     p_{11}>p_{22}+l     \qmbox{or}
  p_{11}< p_{12}-l
  \;,\nonumber\\
  && e_1 \ket{p_{12}}{p_{22}}{p_{11}}' 
  = \Big([p_{12}-p_{11}][p_{11}-p_{22}]\Big)^{1/2}
  \ket{p_{12}}{p_{22}}{p_{11}+1}'  
  \nonumber\\
  &&  \qquad\qmbox{for}  p_{22}+l> p_{11}> p_{12}-l
  \;,\nonumber \\
  && e_1 \ket{p_{12}}{p_{22}}{p_{12}-l} =
  \Big([p_{12}-p_{22}-l]\Big)^{1/2}
  \ket{p_{12}}{p_{22}}{p_{12}-l+1}' \qquad \mbox{(instead of $0$)}
  \;,\nonumber\\
  && e_1 \ket{p_{12}}{p_{22}}{p_{22}+l}' = 0
  \;,\nonumber\\
  && e_1 \ket{p_{22}+l}{p_{12}-l}{p_{22}+l}' = 
  \Big(-[p_{12}-p_{22}-l]\Big)^{1/2}
  \ket{p_{12}}{p_{22}}{p_{22}+l+1} 
  \nonumber\\
  &&\qquad\qquad \mbox{(instead of $0$)}
  \;,\nonumber\\
  && e_1 \ket{p_{22}+l}{p_{12}-l}{p_{11}}' = 
  \Big([p_{12}-p_{11}][p_{11}-p_{22}]\Big)^{1/2}
  \ket{p_{22}+l}{p_{12}-l}{p_{11}+1}' 
  \nonumber\\
  &&\qquad \qquad
  + \Big(-[p_{12}-p_{11}][p_{11}-p_{22}]\Big)^{-1/2} [p_{12}-p_{22}]
  \ket{p_{12}}{p_{22}}{p_{11}+1}  
  \nonumber\\
  &&  \qquad\qmbox{for}  p_{22}+l> p_{11}> p_{12}-l \;.
  \label{eq:indec-e1}
\end{eqnarray}

We can check that $e_1^l$ and $f_1^l$ vanish on this indecomposable
representation.

\section{Action of $e_2$ and $f_2$ in the teepee\label{app:e2f2}}

Let us consider 
$\ket{p_{12}}{p_{22}}{p_{11}}'$ with
$p_{12}-p_{22}>l$, such that both (\ref{eq:triangineq}) and 
(\ref{eq:triangineq2}) are satisfied. 
Then $f_2 \ket{p_{12}}{p_{22}}{p_{11}}'$ if given by a formula
analogous to (\ref{eq:GZ2}), but with primed states on the right hand
side. This is true  as long as the final primed states are defined.
One has
\begin{eqnarray}
  f_2 \ket{p_{22}+l}{p_{12}-l}{p_{11}}' &=&
  \left(\frac {P_1 P_2}{P_3}(1,2;p) 
  \right)^{1/2}
  \ket{p_{22}+l}{p_{12}-l-1}{p_{11}}'
  \nonumber\\  
  & + & 
  \left(\frac {P_1 P_2}{P_3}(2,2;p_{22}) 
  \right)^{1/2}
  \ket{p_{22}+l-1}{p_{12}-l}{p_{11}}'
  \nonumber\\
  & + &  
  \cD_{p_{12}\rightarrow p_{12}-l \atop 
    p_{22}\rightarrow p_{22}+l}
  \left(-\frac {P_1 P_2}{P_3}(1,2;p) 
  \right)^{1/2} 
  \ket{p_{12}-1}{p_{22}}{p_{11}}' 
  \nonumber\\
  & + & 
  \cD_{p_{12}\rightarrow p_{12}-l \atop 
    p_{22}\rightarrow p_{22}+l}
  \left(-\frac {P_1 P_2}{P_3}(2,2;p) 
  \right)^{1/2} 
  \ket{p_{12}}{p_{22}-1}{p_{11}}' 
  \;,
  \label{eq:GZ3}
\end{eqnarray}
\begin{eqnarray}
   e_2  \ket{p_{22}+l}{p_{12}-l}{p_{11}}' & = &
  \left(\frac {P_1 P_2}{P_3}(1,2;p_{12}+1) 
  \right)^{1/2}
  \ket{p_{22}+l}{p_{12}-l+1}{p_{11}}'
  \nonumber\\
  & + & 
  \left(\frac {P_1 P_2}{P_3}(2,2;p_{22}+1) 
  \right)^{1/2}
  \ket{p_{22}+l+1}{p_{12}-l}{p_{11}}'
  \nonumber\\
  & + & 
  \cD_{p_{12}\rightarrow p_{12}-l \atop 
    p_{22}\rightarrow p_{22}+l}
  \left(-\frac {P_1 P_2}{P_3}(1,2;p_{12}+1) 
  \right)^{1/2} 
  \ket{p_{12}+1}{p_{22}}{p_{11}}' 
  \nonumber\\
  & + & 
  \cD_{p_{12}\rightarrow p_{12}-l \atop 
    p_{22}\rightarrow p_{22}+l}
  \left(-\frac {P_1 P_2}{P_3}(2,2;p_{22}+1) 
  \right)^{1/2} 
  \ket{p_{12}}{p_{22}+1}{p_{11}}' 
  \;.
  \label{eq:GZ4}
\end{eqnarray}

where
\begin{equation}
  \cD_{a\rightarrow b} (f) = \lim_{q^l\rightarrow 1}
  \frac1{[l]} \Big(f(a)-f(b) \Big) \;.
  \label{eq:Dlim}
\end{equation}
In the cases we consider, the arguments $a$ and $b$ differ by
multiples of $l$ and the limit in (\ref{eq:Dlim}) is finite. Moreover,
$\cD$ acts as a derivative and 
$\cD_{a\rightarrow b \atop c\rightarrow d}(f) = \cD_{a\rightarrow b}(f)
+ \cD_{c\rightarrow d}(f)$.

\section{
  List of all the possible divergences\label{app:divergences}} 

\subsection{Vanishing denominators in Equation (\ref{eq:GZ2})}

It is easy to see that the denominators in the action of $e_2$ and
$f_2$  (\ref{eq:GZ2}) vanish in
the following cases:
\begin{itemize}
\item For a ``classical'' reason, i.e.\  when $p_{12}-p_{22}=1$ and when
  the action or $e_2$ or $f_2$ would lead to a forbidden state where
  $p_{12}-p_{22}=0$. In such cases, the denominator comes with two
  zeroes in the numerator that cancel this branching.
\item When acting on a \GZ state with $p_{12}-p_{22}=l$. The two 
  resulting states with diverging coefficient actually belong to the
  set of redefined states. The regularization compensates the
  divergence. In the case when one of the resulting states does not
  exist classically, the coefficient of the single remaining state
  (hence not to be redefined) has
  also a zero in the numerator and it remains finite.
\item When the action leads to a state with $p_{12}-p_{22}=l$. 
  If only one initial state can lead to it, the coefficient is finite
  due to a vanishing numerator. If two states lead to it, they have
  the same weight {\em and} same value of $\csl2$, so they are
  redefined. The action on these redefined states contains finite 
  differences of the diverging coefficients.
\end{itemize}

\subsection{Entering and leaving the teepee}

We now summarize the reasons why
the actions of $f_2$ and $e_2$ are well-defined on the boundary of the
teepee. 
Let us first consider the action  of $f_2$ and $e_2$ 
on a state lying  out of the
teepee, the effect of which is to enter the teepee. We have four
different ways of entering:
\begin{itemize}
\item Through the left roof: $p_{22}=p_{13}-l$ for the final
  state. This is reached as $f_2$ acts on
  $\ket{p_{12}}{p_{13}-l+1}{p_{11}}$. A vanishing factor 
  $[p_{13}-p_{22}]$ from the numerator $P_1$
  compensates the diverging factor from the redefinition
  (\ref{eq:transf}) of the final
  state. This is true unless $p_{12}=p_{13}$, in which
  case this vanishing factor compensates  a factor from $P_3$ that
  goes to zero. We arrive in this case 
  on $p_{12}-p_{22}=l$ and there is no  redefinition. 
\item Through the right roof: $p_{12}=p_{33}+l+1$ for the final
  state. This is reached as $e_2$ acts on
  $\ket{p_{33}+l}{p_{22}}{p_{11}}$. A vanishing factor 
  $[p_{12}-p_{33}-1]$ from the numerator $P_1$
  compensates the diverging factor from the redefinition
  (\ref{eq:transf}) of the final
  state. This is true unless $p_{22}=p_{33}+1$, in which
  case this vanishing factor compensates  a factor from $P_3$ that
  goes to zero. We arrive in this case on $p_{12}-p_{22}=l$ and there
  is no  redefinition. 
\item Through the front ``entrance'': $p_{22}=p_{11}-l$ for the final
  state. This is reached as $e_2$ acts on
  $\ket{p_{12}}{p_{11}-l-1}{p_{11}}$, with $p_{12}-p_{11}>0$. 
  A vanishing factor 
  $[p_{11}-p_{22}]$ from the numerator $P_2$
  compensates the diverging factor from the redefinition
  (\ref{eq:transf}) of the final
  state. Note that for $p_{12}=p_{11}$, this final state is not
  redefined, and the compensation comes from the denominator. 
\item Through the back ``entrance'': $p_{12}=p_{11}+l-1$ for the final
  state. This is reached as $f_2$ acts on
  $\ket{p_{11}+l}{p_{12}}{p_{11}}$, with $p_{11}-p_{22}>1$. 
  A vanishing factor 
  $[p_{12}-p_{11}+1]$ from the numerator $P_2$
  compensates the diverging factor from the redefinition
  (\ref{eq:transf}) of the final
  state. Note that for $p_{11}=p_{22}+1$, this final state is not
  redefined, and the compensation comes from the denominator. 
\end{itemize}
 
We now consider $f_2$ and $e_2$ acting on a redefined state, such
that this actions lead to at least one non-redefined state. Again,
the diverging coefficient involved in (\ref{eq:transf}) may be a
source of problem in the boundary. Without entering into details, the
finiteness argument is again that the boundary of the teepee is a
place where one of the numerators $P_1$ or $P_2$ vanishes and
compensates the denominator in (\ref{eq:transf}).


\begin{thebibliography}{99}

\bibitem {Dri} V.G. Drinfeld, 
  {\sl Quantum Groups,} 
  Proc. Int. Congress of Mathematicians, Berkeley, California, Vol.
  {\bf 1}, Academic Press, New York (1986), 798. 
  
\bibitem {Jim} M. Jimbo, {\sl $q$-difference analogue of 
    $\ \cU (\cG )$ and the Yang Baxter equation,} 
  Lett. Math. Phys. {\bf 10} (1985)  63.
  
\bibitem {RosA} M. Rosso, 
  {\sl Finite dimensional representations of the quantum analogue of the 
    enveloping algebra of a complex simple Lie algebra,} 
  Commun. Math. Phys. {\bf 117} (1988) 581.
  
\bibitem {Lusztigi} G. Lusztig,  
  {\sl Quantum deformations of certain simple modules over enveloping
    algebras,} 
  Adv. Math. {\bf 70} (1988) 237.
  
\bibitem {RA}  P. Roche and  D. Arnaudon,
  {\sl Irreducible representations of the quantum analogue of $SU(2)$,}
  Lett. Math. Phys. {\bf 17} (1989) 295.
  
\bibitem {DK}  C. De Concini and V.G. Kac, 
  {\sl Representations of quantum groups at roots of 1,}
  Progress in Math. {\bf 92} (1990) 471 (Birkh{\"a}user).
  
\bibitem {Lusztig} G. Lusztig, 
  {\sl Modular representations and quantum groups,}
  Contemporary Math. {\bf 82} (1989) 59;
  {\sl Quantum groups at roots of $1$,} 
  Geom. Ded. {\bf 35} (1990) 89.
  
\bibitem {DobStAndrews} V.K. Dobrev, 
  in Proc. {\sl Int. Group Theory Conference,} St Andrews, 1989,
  Vol. 1, Campbell and Robertson  (eds.), London 
  Math. Soc. Lect. Notes Series 159, Cambridge University Press, 1991; 
  {\sl Representations of quantum groups,} 
  Proc. ``Symmetries in Science V: Algebraic Structures, their
  Representations, Realizations and Physical Applications'', 
  Schloss Hofen, Austria, (1990), Eds. B.  Gruber, L.C. Biedenharn
  and H.-D. Doebner (Plenum Press, NY, 1991) pp. 93-135.
  
\bibitem {BKW} \v{C}. Burd{\'\i}k, R. C. King and T. A. Welsh, 
  {\sl The explicit construction of irreducible representations of the
    quantum algebras $\cU_q(sl(n))$} 
  Proceedings of the 3rd Wigner Symposium, Oxford, September 1993, 
  Classical and quantum systems, Eds. L. L. Boyle and A. I. Solomon, 
  World Sci., Singapore (1993).
  
\bibitem {DobTruIrreg} V.K. Dobrev and P. Truini,
  {\sl Irregular $U_q(sl(3))$ representations at roots of unity via
    Gelfand--(Weyl)--Zetlin basis,}
  Preprint IC/96/13 and talk given at 
  the XXI International Colloquium on Group Theoretical Methods in
  Physics, 15-20 July 1996, Goslar, Germany. 

\bibitem{JimGZ} M. Jimbo, 
  Lecture notes in Physics 246, Springer (1986), 334.
  
\bibitem {Abdatyp} B. Abdesselam,
  {\sl Special representations of $U_q(sl(N))$ at roots of unity,}
  q-alg/9509016, 
  J. Phys. A: Math. Gen. {\bf 29} (1996) 1201.

\bibitem {ACperiodicsln} D. Arnaudon and A. Chakrabarti,
  {\sl Periodic and partially periodic representations of $SU(N)_{q}$,} 
  Commun. Math. Phys. {\bf 139} (1991) 461.

\bibitem{LasLecThi} A. Lascoux, B. Leclerc and J.-Y. Thibon,
  {\sl Hecke algebras at roots of unity and crystal bases of quantum
  affine algebras,}
  Commun. Math. Phys. {\bf 181} (1996) 205.

\bibitem {PS} V. Pasquier and  H. Saleur, 
  {\sl Common structure between finite systems and conformal field
    theories through quantum groups,} 
  Nucl. Phys. {\bf B330}  (1990) 523.
  
\bibitem {Kel} G. Keller, 
  {\sl Fusion rules of $\ \cU_{q}(SL(2,\CC))$, $q^{m}=1$,} 
  Letters in Math. Phys. {\bf 21} (1991) 273.
  
\bibitem {Suter} R. Suter,
  {\sl Modules over $U_q(sl_2)$,}
  Commun. Math. Phys. {\bf 163} (1994) 359.

\bibitem {Asutroisq} D. Arnaudon, 
  {\sl Periodic and flat irreducible representations of
    $SU(3)_{q}$.}
  Commun. Math. Phys. {\bf 134} (1990) 523.
  
\bibitem {AACuqsln} B. Abdesselam, D. Arnaudon and A. Chakrabarti,
  {\sl Representations of ${\cal U}_q(sl(N))$ at roots of unity,}
  q-alg/9504006,   
  J. Phys. A: Math. Gen. {\bf 28} (1995) 5495.

\bibitem {FGP} P. Furlan, A. C. Ganchev and V. B. Petkova, 
  {\sl Quantum groups and fusion rules multiplicities,} 
  Nucl. Phys. {\bf B343} (1990) 205.
  
\bibitem {Akinetic} D. Arnaudon, {\sl Composition of kinetic momenta:
    the ${\cal U}_q(sl(2))$ case,} 
  hep-th/9212067, 
  Commun. Math. Phys. {\bf 159} (1994) 175.
  
\end{thebibliography}
\end{document}